\documentclass[leqno,12pt]{article}

\usepackage{amssymb}
\usepackage{euscript}
\usepackage[dvips]{graphicx}
\usepackage{flafter}
\usepackage{pstricks}

\usepackage{hyperref}

\usepackage{amsxtra}
 
\usepackage{verbatim}

\usepackage{amsmath}

\usepackage{color}

\usepackage{mathrsfs} 

\usepackage{amsthm}

\evensidemargin\oddsidemargin

\begin{document}

\baselineskip=18pt
\setcounter{page}{1}

\renewcommand{\theequation}{\thesection.\arabic{equation}}
\newtheorem{theorem}{Theorem}[section]
\newtheorem{lemma}[theorem]{Lemma}
\newtheorem{definition}[theorem]{Definition}
\newtheorem{proposition}[theorem]{Proposition}
\newtheorem{corollary}[theorem]{Corollary}
\newtheorem{fact}[theorem]{Fact}
\newtheorem{problem}[theorem]{Problem}
\newtheorem{conjecture}[theorem]{Conjecture}
\newtheorem{claim}[theorem]{Claim}

\theoremstyle{definition} 
\newtheorem{remark}[theorem]{Remark}

\newcommand{\eqnsection}{
\renewcommand{\theequation}{\thesection.\arabic{equation}}
    \makeatletter
    \csname  @addtoreset\endcsname{equation}{section}
    \makeatother}
\eqnsection


\def\r{{\mathbb R}}
\def\e{{\mathbb E}}
\def\p{{\mathbb P}}
\def\z{{\mathbb Z}}
\def\N{{\mathbb N}}
\def\T{{\mathbb T}}

\def\ee{\mathrm{e}}
\def\d{\, \mathrm{d}}



\vglue50pt

\centerline{\Large\bf The maximal drawdown of the Brownian meander}

{
\let\thefootnote\relax\footnotetext{Partly supported by ANR project MEMEMO2 (2010-BLAN-0125).}
}

\bigskip
\bigskip

\centerline{by}

\medskip

\centerline{Yueyun Hu, Zhan Shi and Marc Yor$^\dagger$}

\medskip

\centerline{\it Universit\'e Paris XIII, Universit\'e Paris VI \& Universit\'e Paris VI}

\bigskip
\bigskip
\bigskip

{\leftskip=2.4truecm \rightskip=2.4truecm \baselineskip=15pt \small

\noindent{\slshape\bfseries Summary.} Motivated by evaluating the limiting distribution of randomly biased random walks on trees, we compute the exact value of a negative moment of the maximal drawdown of the standard Brownian meander.

\bigskip


\bigskip


} 

\bigskip
\bigskip

\section{Introduction}
   \label{s:intro}

$\phantom{aob}$Let $(X(t), \, t\in [0, \, 1])$ be a random process. Its maximal drawdown on $[0, \, 1]$ is defined by 
$$
X^\# (1)
:=
\sup_{s\in [0, \, 1]} [\, \overline{X}(s) - X(s)] \, ,
$$

\noindent where $\overline{X}(s) := \sup_{u\in [0, \, s]} X(u)$. There has been some recent research interest on the study of drawdowns from probabilistic point of view (\cite{mijatovic-pistorius}, \cite{nikeghbali}) as well as applications in insurance and finance (\cite{carraro-elkaroui-obloj}, \cite{cheridito-nikeghbali-platen}, \cite{cherny-obloj}, \cite{rieder-wittlinger}, \cite{zhang-hadjiliadis}).

We are interested in the maximal drawdown $\mathfrak{m}^\# (1)$ of the standard Brownian meander $(\mathfrak{m}(t), \, t\in [0, \, 1])$. Our motivation is the presence of the law of $\mathfrak{m}^\# (1)$ in the limiting distribution of randomly biased random walks on supercritical Galton--Watson trees (\cite{yzlocaltree}); in particular, the value of $\e(\frac{1}{\mathfrak{m}^\# (1)})$ is the normalizing constant in the density function of this limiting distribution. The sole aim of the present note is to compute $\e(\frac{1}{\mathfrak{m}^\# (1)})$, which turns out to have a nice numerical value. 

Let us first recall the definition of the Brownian meander. Let $W:= (W(t), \, t\in [0, \, 1])$ be a standard Brownian motion, and let $\mathfrak{g} := \sup\{ t\le 1: \, W(t) =0\}$ be the last passage time at $0$ before time 1. Since $\mathfrak{g}<1$ a.s., we can define
$$
\mathfrak{m}(s)
:=
\frac{|W(\mathfrak{g}+s(1-\mathfrak{g}))|}{(1-\mathfrak{g})^{1/2}} \, ,
\qquad
s\in [0, \, 1] \, .
$$

\noindent The law of $(\mathfrak{m}(s), \, s\in [0, \, 1])$ is called the law of the standard Brownian meander. For an account of general properties of the Brownian meander, see Yen and Yor~\cite{yen-yor}.

\medskip

\begin{theorem}
\label{t:main}

 Let $(\mathfrak{m}(s), \, s\in [0, \, 1])$ be a standard Brownian meander. We have
 \begin{equation}
     \e \Big( \frac{1}{\sup_{s\in [0, \, 1]} [\, \overline{\mathfrak{m}}(s) - \mathfrak{m}(s)]} \Big)
     =
     \Big( \frac{\pi}{2} \Big)^{1/2} \, ,
     \label{cm=}
 \end{equation}
 where $\overline{\mathfrak{m}}(s) := \sup_{u\in [0, \, s]} \mathfrak{m}(u)$.
\end{theorem}

\medskip

The theorem is proved in Section \ref{s:meander}.

We are grateful to an anonymous referee for a careful reading of the manuscript and for many suggestions for improvements.

N.B.\ from the first-named coauthors: This note originates from a question we asked our teacher, {\bf Professor Marc Yor (1949--2014)}, who passed away in January 2014, during the preparation of this note. He provided us, in November 2012, with the essential of the material in Section \ref{s:meander}.

\section{Proof}
\label{s:meander}

$\phantom{aob}$Let $R:= (R(t), \, t\ge 0)$ be a three-dimensional Bessel process with $R(0)=0$, i.e., the Euclidean modulus of a standard three-dimensional Brownian motion. The proof of Theorem \ref{t:main} relies on an absolute continuity relation between $(\mathfrak{m}(s), \, s\in [0, \, 1])$ and $(R(s), \, s\in [0, \, 1])$, recalled as follows.

\medskip

\begin{fact}
\label{f:Imhof}

 {\bf (Imhof~\cite{imhof})}
 Let $(\mathfrak{m}(s), \, s\in [0, \, 1])$ be a standard Brownian meander. Let $(R(s), \, s\in [0, \, 1])$ be a three-dimensional Bessel process with $R(0)=0$. For any measurable and non-negative functional $F$, we have
 $$
 \e \Big[ F(\mathfrak{m}(s), \, s\in [0, \, 1]) \Big]
 =
 \Big( \frac{\pi}{2}\Big)^{\! 1/2} \, 
 \e \Big[ \frac{1}{R(1)} \, F(R(s), \, s\in [0, \, 1]) \Big] \, .
 $$
\end{fact}

\medskip

We now proceed to the proof of Theorem \ref{t:main}. Let 
$$
L
:=
\e \Big( \frac{1}{\sup_{s\in [0, \, 1]} [\, \overline{\mathfrak{m}}(s) - \mathfrak{m}(s)]} \Big) \, .
$$

\noindent Write $\overline{R}(t) := \sup_{u\in [0, \, t]} R(u)$ for $t\ge 0$. By Fact \ref{f:Imhof},
\begin{eqnarray*}
    L
 &=& \Big( \frac{\pi}{2}\Big)^{\! 1/2} \, \e \Big[ \frac{1}{R(1)} \, \frac{1}{\sup_{s\in [0, \, 1]} [\, \overline{R}(s) - R(s)]} \Big] 
    \\
 &=&\Big( \frac{\pi}{2}\Big)^{\! 1/2} \, \int_0^\infty \e \Big[ \frac{1}{R(1)} \, {\bf 1}_{\{ \sup_{s\in [0, \, 1]} [\, \overline{R}(s) - R(s)] < \frac1a \} } \Big] \d a \, ,
\end{eqnarray*}

\noindent the last equality following from the Fubini--Tonelli theorem. By the scaling property, $\e [ \frac{1}{R(1)} \, {\bf 1}_{\{ \sup_{s\in [0, \, 1]} [\, \overline{R}(s) - R(s)] < \frac1a \} } ] = \e [ \frac{a}{R(a^2)} \, {\bf 1}_{\{ \sup_{u\in [0, \, a^2]} [\, \overline{R}(u) - R(u)] < 1 \} } ]$ for all $a>0$. So by means of a change of variables $b=a^2$, we obtain:
$$
L
=
\Big( \frac{\pi}{8}\Big)^{\! 1/2} \, \int_0^\infty \e \Big[ \frac{1}{R(b)} \, {\bf 1}_{\{ \sup_{u\in [0, \, b]} [\, \overline{R}(u) - R(u)] < 1\} } \Big] \d b \, .
$$

\noindent Define, for any random process $X$,
$$
\tau_1^X
:=
\inf\{ t\ge 0: \, \overline{X}(t) - X(t) \ge 1\} \, ,
$$

\noindent with $\overline{X}(t) := \sup_{s\in [0,\, t]} X(s)$. For any $b>0$, the event $\{ \sup_{u\in [0, \, b]} [\, \overline{R}(u) - R(u)] < 1\} $ means $\{ \tau_1^R > b\}$, so
$$
L
=
\Big( \frac{\pi}{8} \Big)^{\! 1/2} \, \int_0^\infty \e \Big[ \frac{1}{R(b)} \, {\bf 1}_{\{ \tau_1^R > b\} } \Big] \d b 
=
\Big( \frac{\pi}{8} \Big)^{\! 1/2} \, \e \Big( \int_0^{\tau_1^R} \frac{1}{R(b)} \d b \Big) ,
$$

\noindent the second identity following from the Fubini--Tonelli theorem. According to a relation between Bessel processes of dimensions three and four (Revuz and Yor~\cite{revuz-yor}, Proposition XI.1.11, applied to the parameters $p=q=2$ and $\nu = \frac12$),
$$
R(t) 
=
U \Big( \frac14 \int_0^t \frac{1}{R(b)} \d b \Big) ,
\qquad
t\ge 0 \, ,
$$

\noindent where $U:= (U(s), \, s\ge 0)$ is a four-dimensional {\it squared} Bessel process with $U(0)=0$; in other words, $U$ is the square of the Euclidean modulus of a standard four-dimensional Brownian motion.

Let us introduce the increasing functional $\sigma(t) := \frac14 \int_0^t \frac{1}{R(b)} \d b$, $t\ge 0$. We have $R=U\circ \sigma$, and
\begin{eqnarray*}
    \tau_1^R
 &=& \inf\{ t\ge 0: \, \overline{R}(t) - R(t) \ge 1\}
    \\
 &=& \inf\{ t\ge 0: \, \overline{U}(\sigma(t)) - U(\sigma(t)) \ge 1\}
    \\
 &=& \inf\{ \sigma^{-1}(s): \, s\ge 0 \hbox{ \rm and } \overline{U}(s) - U(s) \ge 1\} \,
\end{eqnarray*}

\noindent which is $\sigma^{-1}(\tau_1^U)$. So $\tau_1^U = \sigma (\tau_1^R)$, i.e.,
$$
\int_0^{\tau_1^R} \frac{1}{R(b)} \d b = 4 \tau_1^U \, ,
$$

\noindent which implies that
$$
L
=
(2\pi)^{1/2}\, \e(\tau_1^U) \, .
$$

The Laplace transform of $\tau_1^U$ is determined by Lehoczky~\cite{lehoczky}, from which, however, it does not seem obvious to deduce the value of $\e(\tau_1^U)$. Instead of using Lehoczky's result directly, we rather apply his method to compute $\e(\tau_1^U)$. By It\^o's formula, $(U(t)-4t, \, t\ge 0)$ is a continuous martingale, with quadratic variation $4\int_0^t U(s) \d s$; so applying the Dambis--Dubins--Schwarz theorem (Revuz and Yor~\cite{revuz-yor}, Theorem V.1.6) to $(U(t)-4t, \, t\ge 0)$ yields the existence of a standard Brownian motion $B = (B(t), \, t\ge 0)$ such that
$$
U(t)
=
2B(\int_0^t U(s) \d s) + 4 t  \, ,
\qquad
t\ge 0 \, .
$$

\noindent Taking $t:= \tau_1^U$, we get
$$
U(\tau_1^U)
=
2B(\int_0^{\tau_1^U} U(s) \d s) + 4\tau_1^U \, .
$$

\noindent We claim that 
\begin{equation}
    \e\Big[ B(\int_0^{\tau_1^U} U(s) \d s) \Big] 
    =
    0  \, .
    \label{martingale}
\end{equation}

\noindent Then $\e(\tau_1^U) = \frac14 \, \e[U(\tau_1^U)]$; hence
\begin{equation}
    L
    =
    (2\pi)^{1/2}\, \e(\tau_1^U)
    =
    (\frac{\pi}{8})^{1/2}\, \e[U(\tau_1^U)] \, .
    \label{pfeq:1}
\end{equation}

\noindent Let us admit \eqref{martingale} for the moment, and prove the theorem by computing $\e[U(\tau_1^U)]$ using Lehoczky~\cite{lehoczky}'s method; in fact, we determine the law of $U(\tau_1^U)$.

\medskip

\begin{lemma}
\label{l}

 The law of $U(\tau_1^U)$ is given by
 $$
 \p\{ U(\tau_1^U) > a\}
 =
 (a+1) \ee^{-a},
 \qquad \forall a>0.
 $$

\end{lemma}

\medskip

In particular,
$$
\e[U (\tau_1^U)]
=
\int_0^\infty (a+1) \ee^{-a} \d a
=
2.
$$

\noindent Since $L = (\frac{\pi}{8})^{1/2}\, \e[U(\tau_1^U)]$ (see \eqref{pfeq:1}), this yields $L = (\frac{\pi}{2})^{1/2}$ as stated in Theorem \ref{t:main}.

The rest of the note is devoted to the proof of Lemma \ref{l} and \eqref{martingale}.

\bigskip

\noindent {\it Proof of Lemma \ref{l}.} Fix $b>1$. We compute the probability $\p\{ \overline{U} (\tau_1^U) > b\}$ which, due to the equality $\overline{U}(\tau_1^U) = U(\tau_1^U) +1$, coincides with $\p\{ U(\tau_1^U) >b-1 \}$. By applying the strong Markov property at time $\sigma_0^U := \inf\{ t\ge 0: \, U(t)=1\}$, we see that the value of $\p\{ \overline{U} (\tau_1^U) > b\}$ does not change if the squared Bessel process $U$ starts at $U(0)=1$. Indeed, observing that $\sigma_0^U \le \tau_1^U$, $U(\sigma_0^U) =1$ and that $\overline{U}(\tau_1^U) = \sup_{s\in [\sigma_0^U, \, \tau_1^U]} U(s)$, we have 
$$
\p\{ \overline{U} (\tau_1^U) > b\}
=
\p\Big\{ \sup_{s\in [\sigma_0^U, \, \tau_1^U]} U(s) > b\Big\}
=
\p_1 \{ \overline{U} (\tau_1^U) > b\} \, ,
$$

\noindent the subscript $1$ in $\p_1$ indicating the initial value of $U$. More generally, for $x\ge 0$, we write $\p_x (\bullet) := \p( \bullet \, | \, U(0)=x)$; so $\p = \p_0$.

Let $b_0=1 <b_1 < \cdots < b_n := b$ be a subdivision of $[1, \, b]$ such that $\max_{1\le i\le n} (b_i-b_{i-1})\to 0$, $n\to \infty$. Consider the event $\{ \overline{U} (\tau_1^U) > b\}$: since $U(0)=1$, this means $U$ hits position $b$ before time $\tau_1^U$; for all $i\in [1, \, n-1]\cap \z$, starting from position $b_i$, $U$ must hit $b_{i+1}$ before hitting $b_i -1$ (caution: not to be confused with $b_{i-1}$). More precisely, let $\sigma_i^U := \inf\{ t\ge 0: \, U(t) = b_i\}$ and let $U_i(s) := U(s+\sigma_i^U)$, $s\ge 0$; then
$$
\{ \overline{U} (\tau_1^U) > b\}
\subset
\bigcap_{i=1}^{n-1} \{ \hbox{\rm $U_i$ hits $b_{i+1}$ before hitting $b_i -1$}\} \, .
$$

\noindent By the strong Markov property, the events $\{ \hbox{\rm $U_i$ hits $b_{i+1}$ before hitting $b_i -1$}\}$, $1\le i\le n-1$, are independent (caution~: the processes $(U_i(s), \, s\ge 0)$, $1\le i\le n-1$, are not independent). Hence
\begin{equation}
    \p_1\{ \overline{U}(\tau_1^U) > b\}
    \le
    \prod_{i=1}^{n-1} \p_{b_i} \{ \hbox{\rm $U$ hits $b_{i+1}$ before hitting $b_i -1$}\} \, .
    \label{Markov_forte}
\end{equation}

\noindent Conversely, let $\varepsilon>0$, and if $\max_{1\le i\le n} (b_i-b_{i-1})<\varepsilon$, then we also have
$$
\p_1\{ \overline{U}(\tau_{1+\varepsilon}^U) > b\}
\ge
\prod_{i=1}^{n-1} \p_{b_i} \{ \hbox{\rm $U$ hits $b_{i+1}$ before hitting $b_i -1$}\} \, ,
$$

\noindent with $\tau_{1+\varepsilon}^U := \inf\{ t\ge 0: \, \overline{U}(t) - U(t) \ge 1+\varepsilon\}$. By scaling, $\overline{U}(\tau_{1+\varepsilon}^U)$ has the same distribution as $(1+\varepsilon) \overline{U}(\tau_1^U)$. So, as long as $\max_{1\le i\le n} (b_i-b_{i-1})<\varepsilon$, we have
$$
\p_1\{ \overline{U}(\tau_1^U) > b\}
\le
\prod_{i=1}^{n-1} \p_{b_i} \{ \hbox{\rm $U$ hits $b_{i+1}$ before hitting $b_i -1$}\}
\le
\p_1\{ \overline{U}(\tau_1^U) > \frac{b}{1+\varepsilon}\} \, .
$$

\noindent Since $\frac1x$ is a scale function for $U$, we have 
$$
\p_{b_i} \{ \hbox{\rm $U$ hits $b_{i+1}$ before hitting $b_i -1$}\} 
=
\frac{\frac{1}{b_i-1} - \frac{1}{b_i}}{\frac{1}{b_i -1} - \frac{1}{b_{i+1}}} 
=
1 - \frac{\frac{1}{b_i} - \frac{1}{b_{i+1}}}{\frac{1}{b_i -1} - \frac{1}{b_{i+1}}} \, .
$$

\noindent If $\lim_{n\to \infty} \max_{0\le i\le n-1} (b_{i+1}-b_i) = 0$, then for $n\to \infty$,
\begin{eqnarray*}
    \sum_{i=1}^{n-1} \frac{\frac{1}{b_i} - \frac{1}{b_{i+1}}}{\frac{1}{b_i -1} - \frac{1}{b_{i+1}}}
 &=&\sum_{i=1}^{n-1} \frac{b_i-1}{b_i}\, (b_{i+1}-b_i) + o(1)
    \\
 &\to&\int_1^b \frac{r-1}{r} \d r
    \\
 &=&b-1-\log b \, .
\end{eqnarray*}

\noindent Therefore,
$$
\lim_{n\to \infty}
\prod_{i=1}^{n-1} \p_{b_i} \{ \hbox{\rm $U$ hits $b_{i+1}$ before hitting $b_i -1$}\}
=
\ee^{-(b-1-\log b)}
=
b\, \ee^{-(b-1)} \, .
$$

\noindent Consequently,
$$
\p\{ \overline{U}(\tau_1^U) > b\}
=
b\, \ee^{-(b-1)} ,
\qquad \forall b>1.
$$

\noindent We have already noted that $U(\tau_1^U) = \overline{U} (\tau_1^U) -1$. This completes the proof of Lemma \ref{l}.\hfill$\Box$

\bigskip

\noindent {\it Proof of \eqref{martingale}.} The Brownian motion $B$ being the Dambis--Dubins--Schwarz Brownian motion associated with the continuous martingale $(U(t)-4t, \, t\ge 0)$, it is a $(\mathscr{G}_r)_{r\ge 0}$-Brownian motion (Revuz and Yor~\cite{revuz-yor}, Theorem V.1.6), where, for $r\ge 0$, 
$$
\mathscr{G}_r
:=
\mathscr{F}_{C(r)},
\qquad
C(r) 
:=
A^{-1}(r),
\qquad
A(t)
:=
\int_0^t U(s) \d s \, ,
$$

\noindent and $A^{-1}$ denotes the inverse of $A$. [We mention that $\mathscr{F}_{C(r)}$ is well defined because $C(r)$ is an $(\mathscr{F}_t)_{t\ge 0}$-stopping time.] As such,
$$
\int_0^{\tau_1^U} U(s) \d s = A(\tau_1^U) \, .
$$

\noindent For all $r\ge 0$, $\{ A(\tau_1^U) > r\} = \{ \tau_1^U > C(r) \} \in \mathscr{F}_{C(r)} = \mathscr{G}_r$ (observing that $\tau_1^U$ is an $(\mathscr{F}_t)_{t\ge 0}$-stopping time), which means that $A(\tau_1^U)$ is a $(\mathscr{G}_r)_{r\ge 0}$-stopping time. If $A(\tau_1^U)= \int_0^{\tau_1^U} U(s) \d s$ has a finite expectation, then we are entitled to apply the (first) Wald identity to see that $\e[B(A(\tau_1^U))] =0$ as claimed in \eqref{martingale}. 

It remains to prove that $\e[A(\tau_1^U)]<\infty$.

Recall that $U$ is the square of the Euclidean modulus of an $\r^4$-valued Brownian motion. By considering only the first coordinate of this Brownian motion, say $\beta$, we have 
$$
\p\Big\{ \sup_{s\in [0, \, a]} U(s) < a^{1-\varepsilon} \Big\}
\le
\p\Big\{ \sup_{s\in [0, \, a]} |\beta(s)| < a^{(1-\varepsilon)/2} \Big\}
=
\p\Big\{ \sup_{s\in [0, \, 1]} |\beta(s)| < a^{-\varepsilon/2} \Big\} \, ;
$$

\noindent so by the small ball probability for Brownian motion, we obtain:
$$
\p\Big\{ \sup_{s\in [0, \, a]} U(s) < a^{1-\varepsilon} \Big\}
\le 
\exp(- c_1 \, a^\varepsilon) \, ,
$$

\noindent for all $a\ge 1$ et all $\varepsilon\in (0, \, 1)$, with some constant $c_1 = c_1 (\varepsilon) >0$. On the event $\{ \sup_{s\in [0, \, a]} U(s) \ge a^{1-\varepsilon} \}$, if $\tau_1^U > a$, then for all $i\in [1, \, a^{1-\varepsilon} -1] \cap \z$, the squared Bessel process $U$, starting from $i$, must first hit position $i+1$ before hitting $i-1$ (which, for each $i$, can be realized with probability $\le 1-c_2$, where $c_2 \in (0, \, 1)$ is a constant that does not depend on $i$, nor on $a$). Accordingly,\footnote{This is the special case $b_i := i$ of the argument we have used to obtain \eqref{Markov_forte}.}
$$
\p\Big\{ \sup_{s\in [0, \, a]} U(s) \ge a^{1-\varepsilon}, \; \tau_1^U > a \Big\} 
\le
(1-c_2)^{\lfloor a^{1-\varepsilon} -1\rfloor}
\le
\exp(- c_3 \, a^{1-\varepsilon}) \, ,
$$

\noindent with some constant $c_3>0$, uniformly in $a\ge 2$. We have thus proved that for all $a\ge 2$ and all $\varepsilon\in (0, \, 1)$,
$$
\p\{ \tau_1^U > a\} 
\le
\exp(- c_3 \, a^{1-\varepsilon}) + \exp(- c_1 \, a^\varepsilon) .
$$

\noindent Taking $\varepsilon := \frac12$, we see that there exists a constant $c_4>0$ such that
$$ 
\p\{ \tau_1^U > a\} 
\le
\exp(- c_4 \, a^{1/2}) ,
\qquad \forall a\ge 2.
$$

\noindent On the other hand, $U$ being a squared Bessel process, we have, for all $a>0$ and all $b\ge a^2$, 
$$
\p\{ A(a) \ge b\} 
=
\p\{ A(1) \ge \frac{b}{a^2} \} 
\le
\p\Big\{ \sup_{s\in [0, \, 1]} U(s) \ge \frac{b}{a^2} \Big\} 
\le \ee^{- c_5 \, b/a^2} \, ,
$$

\noindent for some constant $c_5>0$. Hence, for $b\ge a^2$ and $a\ge 2$,
$$
\p\{ A( \tau_1^U ) \ge b\}
\le
\p\{ \tau_1^U > a\} + \p\{ A(a) \ge b\}
\le
\exp(- c_4 \, a^{1/2}) + \ee^{- c_5 \, b/a^2} \, .
$$

\noindent Taking $a:= b^{2/5}$ gives that
$$
\p\{ A( \tau_1^U ) \ge b\} \le \exp(- c_6 \, b^{1/5}) \, ,
$$

\noindent for some constant $c_6 >0$ and all $b\ge 4$. In particular, $\e[A( \tau_1^U )]<\infty$ as desired.\hfill$\Box$

\bigskip
\bigskip


{\footnotesize

\baselineskip=12pt

\noindent
\begin{tabular}{lllll}
& \hskip10pt Yueyun Hu
    & \hskip15pt Zhan Shi 
    & \hskip15pt Marc Yor$^\dagger$\\
& \hskip10pt D\'epartement de Math\'ematiques
    & \hskip15pt LPMA, Case 188 
    & \hskip15pt LPMA, Case 188 \& IUF \\
& \hskip10pt Universit\'e Paris XIII
    & \hskip15pt Universit\'e Paris VI
    & \hskip15pt Universit\'e Paris VI \\
& \hskip10pt 99 avenue J-B Cl\'ement
    & \hskip15pt 4 place Jussieu
    & \hskip15pt 4 place Jussieu \\
& \hskip10pt F-93430 Villetaneuse
    & \hskip15pt F-75252 Paris Cedex 05
    & \hskip15pt F-75252 Paris Cedex 05 \\
& \hskip10pt France
    & \hskip15pt France
    & \hskip15pt France \\
& \hskip10pt {\tt yueyun@math.univ-paris13.fr}
    & \hskip15pt {\tt zhan.shi@upmc.fr} 
    & \hskip15pt {\tt my19492014@gmail.com}
\end{tabular}

}

\end{document}